%% file: jam-rev.tex
\documentclass[11pt]{article}

\usepackage{authblk}
\usepackage{amssymb,amsmath,amsthm}
\usepackage{enumerate}
\usepackage{xcolor}

\usepackage{epigraph}

\usepackage{etoolbox}

\makeatletter
\patchcmd{\epigraph}{\@epitext{#1}}{\itshape\@epitext{#1}}{}{}
\makeatother

\newtheorem{theorem}{Theorem}
\newtheorem{lemma}{Lemma}
\newtheorem{corollary}{Corollary}
\newtheoremstyle{neosn}{0.5\topsep}{0.5\topsep}{\rm}{}{\sc}{.}{ }{\thmname{#1}\thmnumber{ #2}\thmnote{ {\mdseries#3}}}
\theoremstyle{neosn}
\newtheorem{proposition}{Proposition}

\newtheorem{definition}{Definition}
\newtheorem{problem}{Problem}

\newcommand{\FOL}{\mathrm{FOL}}
\newcommand{\Sym}{\textit{Sym}}

\newcommand{\WU}{\mathrm{WU}}
\newcommand{\OWU}{\mathrm{OWU}}
\newcommand{\PF}{\mathrm{PF}}
\newcommand{\ignore}[1]{}

\begin{document}

\title{Axiomatizing origami planes}
\author[1]{L. Beklemishev}
\author[2]{A. Dmitrieva}
\author[3]{J.A. Makowsky}
\affil[1]{Steklov Mathematical Institute of RAS, Moscow, Russia}
\affil[1]{National Research University Higher School of Economics, Moscow}
\affil[2]{University of Amsterdam, Amsterdam, The Netherlands}
\affil[3]{Technion -- Israel Institute of Technology, Haifa, Israel}
\date{\today}
\maketitle
\input{jam-abstract}

\input{intro-rev}

\input{wu-rev}
\input{ptrings-rev}

\input{hjaxioms-rev}

\input{ordered-rev}
\input{vieta-rev}

\input{jam-dconclu}
\bibliographystyle{alpha}
\bibliography{jam-origami}
\end{document}

%% file: jam-abstract.tex
\begin{abstract}
We provide a variant of an axiomatization of elementary geometry based on logical axioms in the spirit of Huzita--Justin axioms for the origami constructions. We isolate the fragments corresponding to natural classes of origami constructions such as Pythagorean, Euclidean, and full origami constructions. The set of origami constructible points for each of the classes of constructions provides the minimal model of the corresponding set of logical axioms.

Our axiomatizations are based on Wu's axioms for orthogonal geometry and some modifications of Huzita--Justin axioms. We work out bi-interpretations between these logical theories and theories of fields as described in
J.A. Makowsky (2018). Using a theorem of M. Ziegler (1982)
which implies that the first order theory of Vieta fields is undecidable, we conclude that the first order theory of our
axiomatization of origami is also undecidable.
\end{abstract}

%% file: intro-rev.tex
\epigraph{Dedicated to Professor Dick de Jongh on the occasion of his 81st birthday}

\section{Introduction}
\label{se:intro}


The ancient art of paper folding, known in Japan and all over the world as  \emph{origami}, has a sufficiently long tradition in mathematics.\footnote{The book by M. Friedman \cite{Fri-Orig} is an excellent source on the history of mathematical origami.}  The pioneering book by T. Sundara Rao \cite{Rao} on the science of paper folding attracted attention by Felix Klein. Adolf Hurwitz dedicated a few pages of his diaries to paper folding constructions such as the construction of the golden ratio and of the regular pentagon, \ignore{Among other things he developed an (approximate) construction by paper folding (recorded in his notebooks) of a regular pentagon} see \cite{Fri-Orig}. Early books on origami science, such as Young and Young \cite{YY}, considered paper folding mainly as recreational mathematics and as a means of introducing geometry to children. 

Mathematical origami has been advanced by Margherita Piazzola Beloch in the 1930s. Even though her work was published in Italian at the time of Mussolini and remained essentially unnoticed,  she found several key ideas that were rediscovered only much later, see \cite{Beloch}. In particular, she introduced a new operation, the so-called \emph{Beloch fold}, using which she showed a construction of a segment of length $\sqrt[3]{2}$ and more generally showed how to construct the roots of equations of degrees $3$ and $4$. 

The interest towards mathematical origami reemerged in the 1980s through the work of several enthusiasts such as Humiaki Huzita, who organized the successful \emph{First international meeting on origami science and technology} \cite{Proc-Huzita}. Among various mathematical and computational questions studied in relation to origami in modern literature are whether a given crease pattern can be flat-folded, computational complexity of folding problems, and many others (see e.g.\ \cite{bk:DemRou}).

The general problem of origami design is, given a three-dimensional shape, to find a sequence of folds to create an origami approximating that shape (if possible) from a square sheet of paper. Great practical advances in this problem have been achieved by Robert J.~Lang who developed an algorithm and a program that helped to create folding  patterns for certain basic shapes. These basic shapes are subsequently relatively easy to refine into the required artistic images. Using his algorithms Lang designed origami models of very complicated and realistic nature \cite{LangArt}. 

Apart from the visual art, origami science has impressive practical applications in the design of unfoldable structures, for example, satellite solar batteries, telescope mirrors, and even coronary bypasses.

\bigskip
One of the basic mathematical questions about origami, highly relevant for this paper, is origami constructibility: which lengths can be constructed by foldings from a square sheet of paper or, more generally, from a given set of initial points on the plane? This problem is similar to the classical one about compass and ruler constructions and has been studied from the very beginnings of mathematical origami. To systematically investigate this problem and to prove results on non-constructibility, there was a need to formulate basic rules of the game: a finite set of operations to which one can reduce any complex folding.    

 J. Justin \cite{justin1989resolution} and H. Huzita \cite{Huzita} identified a list of six  such operations (H-1), $\ldots$ , (H-6), later called the Huzita -- Justin or Huzita -- Hatori axioms. (The seventh operation proposed by Justin and Hatori was later shown to be reducible to the other ones.) 
\begin{description}
\item[(H-1):]
Given two points $P_1$ and $P_2$, one can make a fold that passes through both of them.
\item[(H-2):]
Given two points $P_1$ and $P_2$, one can make a fold that places $P_1$ onto $P_2$.
\item[(H-3):]
Given two lines $\ell_1$ and $\ell_2$, one can make a fold that places $\ell_1$ onto $\ell_2$.
\item[(H-4):]
Given a point $P$ and a line $\ell$, one can make a fold orthogonal to $\ell$ that passes through $P$.
\item[(H-5):]
Given two points $P_1$ and $P_2$ and a line $\ell_1$, one can make a fold that places $P_1$ onto $\ell_1$ and passes through $P_2$.
\item[(H-6):]
Given two points $P_1$ and $P_2$ and two lines $\ell_1$ and $\ell_2$, one can make a fold that places $P_1$ onto $\ell_1$ and $P_2$ onto $\ell_2$.
\item[(H-7):] Given a point $P$ and two lines $\ell_1,\ell_2$ one can make a fold orthogonal to $\ell_1$ that places $P$ onto $\ell_2$. 
\end{description}

We note that (H-6) is essentially the Beloch fold mentioned earlier. Some of the expositions also take the points in (H-1) and (H-2) to be distinct and then demand the lines obtained in (H-1), (H-2) and (H-4) to be unique. This obviously does not affect origami constructibility. Therefore, for the reasons of simplicity, we describe the rules without the requirements of uniqueness. 

Lang and Alperin~\cite{AlpLang} showed that (H-1) -- (H-7) can be characterized as all the operations that define a unique fold by alignment of points and finite line segments. Using these operations elegant solutions of the two classical problems --- the doubling of the cube and the angle trisection --- were found.

Alperin \cite{alperin2000mathematical} characterized the classes of points (corresponding to certain subfields of $\mathbb{C}$) constructible using foldings defined by natural subsets of Huzita--Justin operations. Our paper uses the ideas of Alperin in the part that deals with the interpretations between geometric theories and the theories of the corresponding classes of fields.

\bigskip 
As mentioned above, Huzita -- Justin axioms were not meant to be understood as axioms in a logical sense but rather as specifying a (not necessarily deterministic or always defined) set of operations generating the origami constructible points on the plane. The distinction between logical axioms and the operational approach was clear in the early XXth century works on origami mathematics (see Chapter 5 of \cite{Fri-Orig}, especially the discussion on page~285). In fact, it was considered a weakness of the origami approach compared to axiomatic geometry in the tradition of Hilbert. The aim of this paper is to connect the two views and to adapt the principles of origami constructions to serve as logical axioms of planar geometry.

In \cite{makowsky2018undecidability,makowsky2019can}
a proof was outlined of the statement that the first order theory of origami planes was undecidable. Although the proof strategy is feasible, the exact definition of the first order theory of origami planes was left imprecise. In particular, the role of the betweenness relation was overlooked. The purpose of this paper is to provide a precise definition of the first order theory of origami planes and to establish its properties. 

The second aim of the paper is to work out mutual first order interpretations between our logical theories of origami and certain classes of fields, as suggested in \cite{makowsky2019can}. Ziegler's theorem then shows the algorithmic undecidability of the first order theories of these structures. 

Even though these results are not really surprising, there are many details to take care of in their accurate proofs. The proofs rely on certain background from several areas: The development of axiomatic orthogonal geometry, the construction of coordinatization in the context of sufficiently weak axioms of geometry, and the logical techniques of first order interpretations. Since we deal with the questions of first order definability and interpretability, we need to be careful about the details of the constructions usually presented in a less formal way. Therefore, our paper has a long introductory part where we present these background results. We assume that the reader has some basic knowledge of elementary algebra, analytic geometry and first order logic, including the notion of interpretation.

\ignore{
As it turns out, if one directly translates Huzita -- Justin axioms into a logical language
treating the free variables in these axioms as universally quantified, the axioms (H-5) and (H-6) become inconsistent.
Therefore, the question arises, how we can adapt the principles of origami constructions
to serve as logical axioms of planar geometry. There are two issues here: what are the basic predicates,
in which the the axioms (H-1), $\ldots$ , (H-6) can be formulated, and what are the additional
axioms which govern the basic predicates.}

\ignore{In the next section we give background and outline our main results.
The axioms mentioned in this section are given explicitly in Section
\ref{se:pythagoras} for Wu planes, in Section \ref{se:hjaxioms} for the first four Huzita -- Justin axioms, and in Section \ref{se:ordered} for the betweenness axioms.}

The plan of the paper is as follows. In Section \ref{se:wuplanes} we present background on orthogonal geometry using the book by Wu \cite{bk:Wu1994} as our main reference source. We introduce the notion of metric Wu plane which is basic for all further developments.  In Section \ref{interp} we present background on interpretations. In Section \ref{se:ptrings} we present the details of the coordinatization of Wu planes. This is an expanded version of the corresponding section of \cite{makowsky2019can}. In Section \ref{se:pyth} we prove the result on the bi-interpretability between the classes of metric Wu planes and of Pythagorean fields. In Section \ref{se:undec} we state a general undecidability result for geometric theories of metric Wu planes based on Ziegler's theorem.

In Section \ref{se:hjaxioms} we introduce appropriate logical versions of the Huzita -- Justin axioms. We also show that metric Wu planes satisfy the first four Huzita -- Justin axioms. In Section \ref{se:ordered} we discuss the role of order axioms, introduce ordered metric Wu planes and their relation to ordered Pythagorean fields.
In Section \ref{se:vieta} we introduce Euclidean and Vieta fields and prove our main Theorems \ref{th:pyth} and \ref{th:vieta}. Finally, in Section \ref{se:conclu} we discuss the remaining open questions.

%% file: wu-rev.tex
\section{Background on axiomatic geometry}
\label{se:wuplanes}
\label{se:pythagoras}

In this section we present necessary background on axiomatic geometry and begin our discussion of the question 
how one can adapt the principles of origami constructions to serve as logical axioms of planar geometry. There are two issues here: what are the basic predicates, in which the logical analogs of the axioms (H-1), $\ldots$ , (H-6) can be formulated, and what are the additional axioms which govern the basic predicates. 

Looking at the usual statements of (H-1), $\ldots$ , (H-6) we see that these principles use the primitive notions of points, lines (folds), incidence (\emph{``a fold $\ell$ passes through a point $P$''}), which are basic for all standard axiomatizations of elementary geometry. In addition, the predicate \emph{``a fold along $\ell$ identifies points $P$ and $Q$''} and possibly some others are used. 

The notion of a fold along a given line is similar to that of reflection with respect to a line. There are axiomatizations of geometry based on reflection as a basic notion (see \cite{bk:Bach1973}). However, in this paper we find it convenient to rely on the well-developed axiomatizations of geometry based on the notion of orthogonality of lines.

The predicate $\ell\perp m$ \emph{``lines $\ell$ and $m$ are orthogonal''} is sufficiently natural from the point of view of origami constructions: orthogonality can be tested by folding a paper along one of the lines. Moreover, we will see later that the majority of principles corresponding to Huzita--Justin axioms are easily expressible using incidence and orthogonality. Thus, we take Wu's axiomatization of orthogonal plane geometry as basic \cite{bk:Wu1994}. We use this book as a reference source for careful proofs of many technical statements that we will need. 

As we will see, for some of the logical versions of Huzita--Justin axioms, namely (H-3), (H-5) and (H-6), the notion of \emph{betweenness} also plays a role. We are going to analyze the situation in Section \ref{se:ordered}. The predicate $Be(P_1, P_2, P_3)$ \emph{``$P_2$ is between $P_1$ and $P_3$''} holds if $P_2$ belongs to the line segment $P_1P_3$ and all three points are distinct. In order to state the full axiomatization of origami geometry the predicate of betweenness will be added to our vocabulary.


Our axiomatization uses
two basic sorts of variables, denoting lines and points, respectively.
We use $P, P_1, \ldots, P_i$, and more liberally, upper case letters, to denote points and
$\ell, \ell_1, \ldots, \ell_i$, and more liberally, lower case letters, to denote lines.

We use the following basic relations:
\begin{enumerate}[(i)]
\item
the \emph{incidence} relation
$ P \in \ell$ between points and lines,
\item
the \emph{orthogonality} relation $\ell_1 \perp \ell_2$ between two lines,
\item
and the \emph{betweenness} (or \emph{order}) relation between three points $P_1, P_2, P_3$ denoted by
$Be(P_1, P_2, P_3)$.
\end{enumerate}

The two-sorted language can be considered as a notational variant of a single-sorted one, as
the set of lines and the set of points are definable using the incidence relation:
$\ell$ is a line iff $\exists P (P \in \ell)$,
and
$P$ is a point iff $\exists \ell (P \in \ell)$. This allows us to apply various standard notions and results for (one-sorted) first order logic in our context. 

\ignore{We also note that the equidistance relation $PQ\cong RS$ between points $P,Q$
and $R,S$ is definable from incidence and orthogonality~\cite[page 25]{bk:Wu1994}.}

We denote
by $\tau_{wu}$ the vocabulary corresponding to (i)-(ii)
and
by $\tau_{o-wu}$ the vocabulary corresponding to (i)-(iii). By
$\FOL(\tau_{wu})$
and
$\FOL(\tau_{o-wu})$ we denote
the corresponding sets of first order formulas.

Next we turn to the axioms of two-dimensional orthogonal geometry as presented in \cite{bk:Wu1994}. The axioms are subdivided into several groups. 

\paragraph{Hilbert's axioms of incidence}
\begin{description}
\item[(I-1):]
For any two distinct points $A,B$ there is a unique line $\ell$
with $A \in \ell$ and $B\in \ell$.
\item[(I-2):]
Every line contains at least two distinct points.
\item[(I-3):]
There exist three distinct points $A, B, C$ such that no line $\ell$
contains all of them.
\end{description}

\paragraph{Hilbert's (sharper) axiom of parallels}
\begin{description}
\item[(ParAx):]
Let $\ell$ be any line and $A$ a point not on $\ell$. Then there exists one and only
one line determined by $\ell$ and $A$ that passes through $A$ and does not intersect $\ell$.
\end{description}

\paragraph{Axiom  schema of infinity and Desargues' axioms}
\begin{description}
\item[(InfLines):]
Given  distinct non-collinear $A, B, C$ and $\ell$ with $A \in \ell$, $B, C \not\in \ell$ we construct a line $\ell_1$ going through $C$ and parallel to $AB$ and define $A_1$ as the intersection of  $\ell_1$ and $\ell$. Inductively, we define $\ell_n$ as a line going through $C$ and parallel to $A_nB$ and define $A_{n+1}$ as its intersection with $\ell$.
Then all the  $A_i$ are distinct.

\item[(De-1):]
If the three pairs of the corresponding sides of two
triangles $ABC$ and $A' B'C'$ are all parallel to each other, i.e.,
$AB \parallel A'B'$, $AC \parallel A'C'$, $BC \parallel B'C'$,
then the three lines $AA'$, $BB'$, $CC'$ joining the corresponding vertices of these
two triangles are either parallel to each other or concurrent.
\item[(De-2):]
If two pairs of the corresponding sides of two triangles
$ABC$ and $A' B'C'$ are parallel to each other, say
$AB \parallel A'B'$, $AC \parallel A'C'$,
and the three lines joining the corresponding vertices are distinct yet either
concurrent or parallel to each other, then the third pair of the corresponding
sides are also parallel to each other, i.e.,
$BC \parallel B'C'$.
\end{description}

\begin{definition}
A $\tau_{\in}$ structure $\Pi$ is a \emph{Desarguesian plane} if it satisfies (I-1, I-2, I-3), the axiom of infinity (InfLines), (ParAx) and the two axioms of Desargues (De-1) and (De-2).
\end{definition}

In order to introduce the orthogonality axioms in the plane, we consider a new relation of orthogonality $\perp$ and the language $\tau_{wu}$ consisting of $\in$ and $\perp$.

\paragraph{Orthogonality  axioms}
\begin{description}
\item[(O-1):]
$\ell_1 \perp \ell_2$ iff $\ell_2 \perp \ell_1$.
\item[(O-2):]
Given $O$ and $\ell_1$, there exists exactly one line $\ell_2$ with
$\ell_1 \perp \ell_2$ and $O \in \ell_2$.
\item[(O-3):] If $\ell_1 \perp \ell_2$ and $\ell_1 \perp \ell_3$
then $\ell_2  \parallel \ell_3$ or $\ell_2 = \ell_3$.
\item[(O-4):]
For every $O$ there is an $\ell$ with
$O \in \ell$ and $\ell \not \perp \ell$.
\item[(O-5):]
The three heights of a triangle intersect in one point.
\end{description}

Concerning Axiom (O-4) we remark that lines $\ell$ such that $\ell\perp \ell$ are called \emph{isotropic}.

{\bf Caveat:}
Without any axioms of order, the axioms of metric Wu planes do not exclude the existence of isotropic lines. Since our analysis concerns the role of the axioms of order in the statements of origami principles, we do not want to assume that there are no isotropic lines outright. 

\begin{definition}
A $\tau_{wu}$ structure $\Pi$ is an \emph{orthogonal Wu plane} if it is a Desarguesian plane satisfying orthogonality axioms (O-1, O-2, O-3, O-4, O-5).
\end{definition}

\paragraph{Axiom of symmetric axis.}
We assume that we are working in an orthogonal Wu plane. In order to formulate the next axiom, we follow \cite[page 22, Definition 3]{bk:Wu1994} to define the relation of being a symmetric point using only the Incidence relation.

For two arbitrary points $A \neq B$ on a line $\ell$, take an arbitrary $E \notin \ell$ and construct a line $\ell'$ parallel to $\ell$ such that $E \in \ell'$. Let $D$ be the intersection of $\ell'$ and a line going through $B$ parallel to $AE$. Then $ABDE$ is a parallelogram. Finally, construct $C$ as the intersection of $\ell$ and the line going through $D$ and parallel to $EB$. Then, due to Desargues' axioms, point $C$ is independent of
the choice of $E$. We say that $C$ is the \emph{symmetric point of $A$ with respect to $B$}. In addition, for any point $P$ we say that $P$ is a symmetric point of $P$ with respect to $P$. 


Next we define a notion of a midpoint, following \cite[page 23, Definition 4]{bk:Wu1994}. Let $A, B$ be two points on a line $\ell$. Draw through $A$ a line $\ell'$ distinct from $\ell$ and take thereon a point $M'$ distinct from $A$. Construct the symmetric point $B'$ of $A$ with respect to $M'$ and draw through $M'$ a line $M'M \parallel B'B$, meeting $\ell$ at $M$. Then $M$ is independent of the choice of $\ell'$, $M'$, and is called the \emph{midpoint} of $A$ and $B$. Define the midpoint of two coincident points to be the point itself.

Now we use the relation of orthogonality to define symmetric axis following \cite[page 75]{bk:Wu1994}. For any pair $A, B$ of two distinct
points, let the unique line through the midpoint of $A$ and $B$ and perpendicular to the
line $AB$ be the \emph{perpendicular bisector of $A, B$}. Clearly, if $AB$ is an isotropic
line, then its perpendicular bisector is $AB$ itself.

Let the perpendicular bisector of $A, B$ be $\ell$. We
call $A$ the \emph{symmetric point of $B$ with respect to $\ell$} or \emph{$\ell$ the symmetric axis of
$A, B$}. Any point $A$ on $\ell$ is said to be a symmetric point of itself with respect
to $\ell$. (Whenever $\ell$ is isotropic, it is the symmetric axis for any two points $A, B$ on $\ell$.) We denote this relation as $Sym(A, \ell, B)$. It will be important in our treatment of origami geometry.

Let $\ell_1$ be any line and $\ell$ be a non-isotropic line.
By \cite[page 76, Property 1]{bk:Wu1994}, the points symmetric to the points from $\ell_1$ with respect to $\ell$ are
also lying on a unique line, say $\ell_2$. In this case we call $\ell_2$ the
\emph{symmetric line of $\ell_1$ with respect to $\ell$}, or $\ell$ a \emph{symmetric axis of $\ell_1$ and $\ell_2$}.

Finally, we can add the following axiom to the list.



\begin{description}
\item[(AxSymAx):]
Any two intersecting non-isotropic lines have a symmetric axis.
\end{description}

Now we are ready to give

\begin{definition}
A $\tau_{wu}$ structure $\Pi$ is a \emph{metric Wu plane} if it is an orthogonal Wu plane satisfying (AxSymAx).
\end{definition}

Following Wu, we use the word `metric' here, since assuming the axiom of symmetric axes gives us some metric properties, in particular the notion of congruence.

%% file: ptrings-rev.tex
\section{Background on interpretations} \label{interp}

We assume familiarity of the reader with the notion of first-order interpretation, we follow the terminology and notations of Makowsky~\cite{makowsky2019can} here. More background on interpretations can be found, for example, in \cite{Sho67,bk:Hodges93,FrVi,ar:MakowskyTARSKI}. Particular interpretations that we consider here are many-dimensional relative interpretations with parameters. 
In this section we briefly remind the notions of translation scheme, interpretation and bi-interpretability of first order theories. The readers already familiar with these concepts can skip it and move on to the next section.

\begin{definition}
Let $K$ and $L$ be first-order relational signatures, $n$ a positive integer, $\vec p$ a list of variables. An ($n$-dimensional) \emph{translation scheme} $\Sigma$ with definable parameters $\vec{p}$ from $L$ to $K$ is specified by three items:
\begin{itemize}
    \item[-] a formula $\delta_\Sigma(x_0, \ldots, x_{n-1}, \vec{p})$ of signature $K$,
    \item[-] for each relational symbol $R (y_0, \ldots, y_{m-1})$ of $L$, a formula $R_\Sigma(\vec{x}_0, \ldots, \vec{x}_{m-1}, \vec{p})$ of $K$ in which the $\vec{x}_i$ are disjoint $n$-tuples of distinct variables,
    \item[-] a formula $\psi_\Sigma(\vec{p})$ in $K$.
\end{itemize}
    Intuitively, $\delta_\Sigma$ denotes the domain of $\Sigma$; $R_\Sigma$ denote the relations of $\Sigma$ and $\psi_\Sigma(\vec{p})$ denotes the admissible range of parameters of $\Sigma$.
    
    In case $L$ is multi-sorted we may use a $\delta_\Sigma$ for each sort. Then the defining formulas for different sorts can also have different dimensions.  
\end{definition}

\ignore{
In our case the translation $RF_{field}$ from the language of fields to $\tau_{wu}$ has a domain formula $x_0 \in \ell_0$, other defining formulas $add_T$, $mult_T$, $x_0 = 0$, $x_0 = 1$ and a parameter defining formula based on our restrictions on the parameters $\ell_0, m_0, Z_0$. The other translation $PP_{wu}$ from $\tau_{wu}$ to the language of fields has two domain formulas: $x_0 = x_0 \ \wedge \ x_1 = x_1$ (or any other tautology in $x_0, x_1$) defining points and $a_0 \neq 0 \vee a_1 \neq 0 \vee a_2 \neq 0$ defining lines. Other defining formulas are
\begin{itemize}
    \item[-] $x_0 = y_0 \wedge x_1 = y_1$ for $x = y$,
    \item[-] $\exists k ( a_0 = k \cdot b_0 \ \wedge \ a_1 = k \cdot b_1 \ \wedge \ a_2 = k \cdot b_2)$ for $\ell = \ell'$,
    \item[-] $a_0 \cdot x_0 + a_1 \cdot x_1 + a_2 = 0$ for $x \in \ell$,
    \item[-] $a_0 \cdot b_0 + a_1 \cdot b_1 = 0$ for $\ell_1 \perp \ell_2$.
\end{itemize}
}

Any translation scheme $\Sigma$ from $L$ to $K$ naturally defines a map from $L$-formulas to $K$-formulas by making it commute with the propositional connectives and relativizing the quantifiers to the domain specified by $\delta_\Sigma$. We denote this map by $\phi \mapsto \phi^\Sigma$.

\begin{definition}
Let $T$ be a $K$-theory and $S$ be an $L$-theory. A translation scheme $\Sigma$ from $L$ to $K$ is an \emph{interpretation of $S$ in $T$} if 
\begin{enumerate}
    \item $T$ proves $\exists \vec p\:\psi(\vec p)$ and $\forall \vec p\:(\psi(\vec p)\to \exists x_0\dots \exists x_{n-1}\delta(x_0,\dots,x_{n-1},\vec p))$.
    \item For any $L$-sentence $\phi$, $$S \vdash \phi \  \Rightarrow \ T \vdash \forall \vec{p} \ (\psi(\vec{p}) \rightarrow \phi^\Sigma).$$ 
\end{enumerate}
\end{definition}

Let $\Sigma$ be an interpretation of $S$ in $T$. Let $M$ be a model of $T$, and let $\vec{m}$ be a tuple of elements of $M$ such that $M \models \psi(\vec{m})$. The interpretation $\Sigma$ naturally defines an internal model $\Sigma^*(M,\vec m)$ of $S$ within $M$. The domain of $\Sigma^*(M,\vec m)$ is the quotient of the set $D:=\{\vec x\in M: M\models \delta(\vec x,\vec m)\}$ by the equivalence relation defined on $D$ by the interpreted equality relation. Formulas $R_\Sigma$ then define the evaluation of the signature of $S$ in $\Sigma^*(M,\vec m)$. In fact, $\Sigma^*(M,\vec m)$ will be  a model of $S$, since for all $\alpha \in S$ we have $T \models \alpha^\Sigma$.

We can often ignore the dependence of  $\Sigma^*(M,\vec m)$ on the choice of parameters $\vec m$ and will denote it by $\Sigma^*(M)$. Next we introduce the notion of bi-interpretability of theories (and of the corresponding classes of models).

\begin{definition}
Let $T$ be a $K$-theory and $S$ be an $L$-theory. Suppose there are interpretations $\Sigma$ of $S$ in $T$ and $\Xi$ of $T$ in $S$. Then, for each model $M$ of $T$, $\Xi^*(\Sigma^*(M))$ is a model of $T$ and, for each model $N$ of $S$, $\Sigma^*(\Xi^*(N))$ is a model of $S$. We demand that models in these pairs are isomorphic and, moreover, definably isomorphic, in the following sense. 

Suppose there are formulas $\alpha$ in $K$ and $\beta$ in $L$, such that for any model $M$ of $T$, $\alpha$ defines in $M$ an isomorphism between $M$ and $\Xi^*(\Sigma^*(M))$, for any admissible choice of parameters, and, for any model $N$ of $S$, $\beta$ defines in $N$ an isomorphism between $B$ and $\Sigma^*(\Xi^*(N))$, for any admissible choice of parameters. Then we say that $T$ and $S$ are \emph{bi-interpretable}.
\end{definition}

Bi-interpretation is a rather strong form of equivalence of theories, in particular it reduces the decision problem for one theory to another.\footnote{There also are well-known weaker conditions under which interpretations preserve decidability of theories. For example, it is sufficient to require that each model of $S$ is isomorphic to a model of the form $\Sigma^*(M)$ for some model $M$ of $T$, see Lemma 4.1 in \cite{makowsky2017can}.} 

\begin{proposition} If $S$ and $T$ are bi-interpretable and $S$ is decidable, then so is $T$. 
\end{proposition}  

We will give concrete examples how this works below. 

\section{Background on coordinatization}
\label{se:ptrings}

\begin{definition}
Pythagorean field is a field for which every sum of two squares is a square:
$$\forall x, y \ \exists z \ x^2 + y^2 = z^2$$
\end{definition}

We would like to describe mutual interpretations between the classes of metric Wu planes and of Pythagorean fields. 

Given a field $\mathcal{F}$ one can define a plane $\Pi$ in a standard manner via Cartesian coordinates, 
we denote such a $\Pi$ as $PP^*_{wu}(\mathcal{F})$, following \cite{makowsky2019can}. It is rather obvious that the two sorts (points and lines) and the predicates of incidence and orthogonality will be definable in the language of fields, which in fact yields a first-order (two-dimensional) interpretation.
Notice that this interpretation is parameter-free.

On the other hand, given a plane $\Pi$ we can 
define a field $\mathcal{F}$ in it by a classical construction known as `coordinatization'. This construction was explicitly introduced by M. Hall in \cite{hall1943projective}.
M. Hall credits \cite{von1857beitrage,hilbert1971foundations} for the original idea.
A good exposition can be found in \cite{blumenthal1980modern,szmielew1983affine}.

Since coordinatization is more complicated than the opposite map, it is not immediately obvious that it yields a first-order interpretation. In fact, the amount of work needed to see it is considerable. The purpose of this section is to give a logically-minded reader enough details, without writing out all the first-order formulas explicitly, so that he or she could be convinced that it, indeed, does.
We follow here \cite{ivanov2016affine},
which contains a particularly nice exposition of this construction.

Let $\Pi$ be a stucture 
satisfying (I-1, I-2, I-3) and (ParAx),
with two distinguished intersecting lines $\ell_0, m_0$ in $\Pi$.
Let $O$ be the point of intersection of the lines $\ell_0$ and $m_0$.
Take any point $Z_0$ such that $Z_0 \notin \ell_0 \cup m_0$.
All the constructions below will depend on a particular choice of parameters $\ell_0$, $m_0$, and $Z_0$ in the above configuration.

\newcommand{\bij}{\textit{bij}}
\begin{lemma}
\label{bijection}
There is a formula $\bij(X,Y,\ell_0, m_0, Z_0) \in \FOL_{\in}$ which,
for every choice of $\ell_0$, $m_0$ and $Z_0$ as above, defines a bijection between the points of $\ell_0$ and  of $m_0$.
\end{lemma}
\begin{proof}
Let $d$ be the line going through the points $O$ and $Z_0$. Let $X \in \ell_0$ and $h(X)$ be the point at the intersection of $d$ and the line $m_X$ parallel to $m_0$ containing $X$.
Let $y(X) \in m_0$ be the point at the intersection of $m_0$ and the line $\ell_X$ parallel to $\ell_0$ and containing $h(X)$.
Clearly $f: \ell_0 \rightarrow m_0$ given by $f(X) =y(X)$ is a bijection and is $\FOL$ definable
by a formula $\bij(X,Y,\ell_0, m_0, Z_0)$.
\end{proof}

We will define a structure $RF^*_{field}(\Pi)$ (this notation is once again taken from \cite{makowsky2019can}) whose universe will be denoted $K$, which we take to be ${\{A: A\in\ell_0\}}$.
Thinking of $\ell_0$ and $m_0$ as the axes of a {\em coordinate system} we can identify the points of $\Pi$
with elements of $K^2$ as follows.

Let $P$ be a point of $\Pi$. The projection of a point $P$ onto $ \ell_0$ is defined by the point $X \in \ell_0$
which is the intersection of the line $m_P$ parallel to $m_0$ with $P \in m_P$.
After analogously projecting a point $P$ onto $m_0$, we get a point $Y \in m_0$. Then the coordinates of $P$ is a pair $(X, f(Y)) \in K^2$.

We define $0$ and $1$ in $K$ by saying that the point $O$ has coordinates $(0,0)$ and the point $Z_0$ has coordinates $(1, 1)$.
Coordinates and elements of $K$ will also be denoted by lower case letters.
It should be clear from the context whether lower case letters denote lines
or elements of $K$.

Next we define the {\em slope} $sl(\ell) \in K \cup \{\infty\}$ of a line $\ell$ in $\Pi$.
If $\ell$ is parallel to $\ell_0$, its slope is $0$ and it is called a {\em horizontal} line.
If $\ell$ is parallel to $m_0$, its slope is $\infty$ and it is called a {\em vertical} line.
For $\ell$ not vertical, let $\ell_1$ be the line parallel to $\ell$ and passing through $0$.
Let $(1,a)$ be the coordinates of the intersection of $\ell_1$ with the vertical line $\ell_2$
passing through $(1,0)$.
Then we define the slope by $sl(\ell)= a \in K$.

This shows:
\begin{lemma}
\label{slope}
There is a first order formula $slope(\ell, A, \ell_0, m_0, Z_0) \in \FOL_{\in}$ which expresses
$sl(\ell)=A$, for any choice of parameters $\ell_0, m_0, Z_0$.
There is also a first order formula $slope_{\infty}(\ell, m_0) \in \FOL_{\in}$ expressing $sl(\ell)=\infty$.
\end{lemma}

\begin{lemma}
\label{slope-1}
\begin{enumerate}[(i)]
\item
Two lines $\ell, \ell_1$ have the same slope, $sl(\ell)=sl(\ell_1)$ iff they are parallel.
\item
For the line $d$, defined in Lemma \ref{bijection}, we have $sl(d)=1$ (because $(1,1) \in d$).
\end{enumerate}
\end{lemma}

We now define a ternary operation $T: K^3 \rightarrow K$ on the set $K = \{A:A\in\ell_0\}$.
We think of $T(a,x,b)= \langle ax+b \rangle$ as the result of multiplying $a$ with $x$ and then adding $b$.
But we yet have to define multiplication and addition.

Let $a,b,x \in K$. Let $\ell$ be the unique line with $sl(\ell)=a \neq \infty$
intersecting the line $m_0$ at the point $P_1$ with coordinates $P_1 =(0,b)$.
Let $\ell_1 =\{ (x,z) \in K^2 : z \in K\}$. The line $\ell$ intersects $\ell_1$
at a unique point, say $P_2= (x,y)$.
We set $T(a,x,b) = y$.

\newcommand{\Ter}{\textit{Ter}}
\begin{lemma}
\label{ptr}
There is a formula $\Ter(a,x,b,y,\ell_0, m_0, Z_0) \in \FOL_{\in}$,
where $a,b,x,y$ range over coordinates and $\ell_0, m_0, Z_0$ are parameters of lines and points,
which expresses that $T(a,x,b)=y$.
\end{lemma}

\begin{lemma}
\label{ptr-1}
The ternary operation $T(a,x,b)$ has the following properties and interpretations:
\begin{description}
\item[(T-1):]
$T(1,x,0)=T(x,1,0)=x$
\\
\emph{$T(1,x,0)=x$ means that the auxiliary line $d =\{ (x,x) \in K^2 : x \in K \}$
is a line with $sl(d)=1$.
\\
$T(x,1,0)=x$ means that the slope
of the line $\ell$ passing through $(0,0)$ and $(1,x)$ is given by $sl(\ell)=x$.
}
\item[(T-2):]
$T(a,0,b)=T(0,a,b)=b$
\\
\emph{The equation $T(a,0,b)=b$ means that the line $\ell$ defined by $T(a,x,b)=y$
intersects $m_0$ at $(0,b)$ (which is the meaning of $ax+b$ in analytic geometry).
\\
The equation $T(0,a,b)=b$ means that the horizontal line $\ell_1$ passing through $(0,b)$
consists of the points $\{ (a,b) \in K^2 : a \in K \}$.}
\item[(T-3):]
For all $a,x,y \in K$ there is a unique $b \in K$ such that $T(a,x,b) =y$
\\
\emph{This means that for every slope  $s$ different from $\infty$ there is a unique line
$\ell$ with $sl(\ell)=s$ passing through  $(x,y)$.}
\item[(T-4):]
For every $a, a', b, b' \in K$ and $a \neq a'$ the equation $T(a,x,b) = T(a'x,b')$ has a unique
solution $x \in K$.
\\
\emph{This means that two lines $\ell_1$ and $\ell_2$ with different slopes not equal to $\infty$
intersect at a unique point $P$.}
\item[(T-5):]
For every $x, y, x', y' \in K$ and $x \neq x'$ there is a unique pair $a,b \in K$ such that
$T(a,x,b)=y$ and
$T(a,x',b)=y'$.
\\
\emph{This means that any two points $P_1, P_2$ not on the same vertical line are contained
in a unique line $\ell$ with
slope different from $\infty$.}
\end{description}
\end{lemma}

A structure $\langle K, T_K, 0, 1 \rangle$  with a ternary operation $T_K$ and $0, 1 \in K$ satisfying (T-1)--(T-5)
is called a {\em planar ternary ring} PTR.
We also define addition $add_T(a,b,c)$ by $T(a, 1, b) =c$
and multiplication $mult_T(a,x,c)$ by $T(a, x, 0) =c$.

Following \cite{makowsky2019can},  we denote the structure $(K; add_T, mult_T, 0, 1)$ as $RF^*_{field}(\Pi)$.
It is shown in \cite{hilbert2013grundlagen}, that if $\Pi$ is a Desarguesian plane,
then $RF^*_{field}(\Pi)$ is a skew-field (a field,
in which the commutativity of the multiplication is not assumed) of characteristic 0.
Moreover, as proved in \cite[page 42, Theorem 1]{bk:Wu1994}, for any such $\Pi$ and any two choices of $\ell_0, m_0, Z_0$,
there is an isomorphism between the two obtained skew-fields.

\section{Metric Wu planes and Pythagorean fields} \label{se:pyth}

The following theorem is stated in a somewhat weaker form in \cite{makowsky2019can},
however it is based on the previous classical results, in particular, by Hall \cite{hall1943projective}
and Wu \cite{bk:Wu1994}.

\begin{theorem} \label{Mak}
\begin{itemize}
\item[(i)] Let $\mathcal{F}$ be a Pythagorean field of characteristic $0$.
Then $PP^*_{wu}(\mathcal{F})$ is a metric Wu plane.
\item[(ii)] Let $\Pi$ be a metric Wu plane. Then $RF^*_{field}(\Pi)$ is a Pythagorean field of characteristic $0$.
\item[(iii)] $RF^*_{field}(PP^*_{wu}(\mathcal{F}))$ is isomorphic to $\mathcal{F}$.
\item[(iv)] $PP^*_{wu}(RF^*_{field}(\Pi))$ is isomorphic to $\Pi$.
\item[(v)] The isomorphisms in (iii) and (iv) are definable, that is, form a bi-interpretation between the classes of metric Wu planes and of Pythagorean fields. 
\end{itemize}
\end{theorem}

\begin{proof}
(i) Axioms (I-1, I-2, I-3) and (ParAx) are shown in~\cite[Proposition 14.1]{bk:Hartshorne2000}.
The infinity axiom holds, since $\mathcal{F}$ has characteristic $0$.
Considering Desargues' axioms, Proposition 14.4 in~\cite{bk:Hartshorne2000} shows that Pappus theorem holds
in a plane defined over a field. Then by Hessenberg's theorem~\cite[page 67]{bk:Wu1994}, Desargues' axioms also hold.

We naturally define lines $a_0x + b_0y + c_0 = 0$ and $a_1x + b_1y + c_1 = 0$ to be orthogonal if $a_0a_1 + b_0b_1 = 0$.
Then  the axiom (O-1) holds by commutativity of multiplication.
The axioms (O-2) and (O-3) hold since we are able to solve systems of linear equations.

Considering axiom (O-5), let $(a_1, b_1)$, $(a_2, b_2)$ and $(a_3, b_3)$ be the coordinates of a triangle. Then a line going through $(a_1, b_1)$ and perpendicular to the opposite side is defined by an equation 
$$ \frac{x - a_1}{b_3 - b_2} + \frac{y - b_1}{a_3 - a_2} = 0.$$
Then since linear equations are solvable in a Pythagorean field of characteristic 0, we can write equations for all three heights and see that the orthocenter $(x,y)$ exists. Therefore, (O-5) holds.

We know that $1 \neq 0$ and hence any line of the form $x = c$ is non-isotropic.
Then for any point $(x_0, y_0)$ there is a non-isotropic line $x = x_0$ passing through it and (O-4) holds.

Suppose an angle is formed by two non-isotropic lines given by $ l_{1}x+m_{1}y+n_{1}=0$ and $l_{2}x+m_{2}y+n_{2}=0$.
Then the internal and external bisectors of the angle are given by the two equations
$$
\frac {l_{1}x+m_{1}y+n_{1}}{\sqrt {l_{1}^{2}+m_{1}^{2}}}=\pm {\frac {l_{2}x+m_{2}y+n_{2}}{\sqrt {l_{2}^{2}+m_{2}^{2}}}}.
$$
Since $\mathcal{F}$ is Pythagorean, the roots exist and are not equal to $0$,
because the lines are non-isotropic. Therefore, bisectors exist and (AxSymAx) holds.

(ii) This statement is extensively discussed in \cite{bk:Wu1994}.
As mentioned above, Hilbert showed in \emph{Grundlagen der Geometrie} that $RF^*_{field}(\Pi)$ forms a
skew-field of characteristic $0$.
Wu first proves in \cite[Section 2.1]{bk:Wu1994} that Linear Pascalian axiom (a version of Pappus theorem)
is sufficient to obtain the commutativity of multiplication and then on \cite[page 72]{bk:Wu1994} shows
that Linear Pascalian axiom holds in any metric Wu plane.

Finally, to conclude that $RF^*_{field}(\Pi)$ is Pythagorean,
we refer to the Pythagorean Theorem (known
as `kou-ku theorem' \cite{bk:Wu1994}) proved on \cite[page 97]{bk:Wu1994}.

(iii) As mentioned above, it is shown in \cite[page 42, Theorem 1]{bk:Wu1994},
that any choice of parameters $\ell_0$, $m_0$ and $Z_0$ in the definition of $RF^*_{field}(\Pi)$ gives us isomorphic fields. Wu provides an explicit construction of such an isomorphism and we shortly describe it here. Let $\ell_0$, $m_0$, $Z_0$ and $\ell'_0$, $m'_0$, $Z'_0$ be two choices of parameters and let $O$, $I$ and $O'$, $I'$ be the corresponding zero and one points. First consider the case when $O = O'$ but $\ell_0 \neq \ell'_0$. Then as shown by Wu the parallel projection with respect to the line $II'$ defines the required isomorphism. Now consider the case $\ell_0 \parallel \ell'_0$ and $OO' \parallel II'$. Then the parallel projection with respect to the line $OO'$ defines the required isomorphism. Finally, observe that any other case can be reduced to those two by taking compositions of isomorphisms. For example, in a general case of $\ell_0$ and $\ell'_0$ being neither coincident nor parallel and their intersection point not equal to $O$ or $O'$, construct a line $\ell''_0$ through $O$ parallel to $\ell'_0$. Then take $I''$ to be a point on $\ell''_0$ such that $O'O \parallel I'I''$. The line $\ell''_0$ and points $O$, $I''$ define a field by taking $m''_0 \perp \ell''_0$, $O \in m''_0$ and $Z''_0$ defined by $I''$. Now using the first case obtain an isomorphism between fields of $\ell_0$ and $\ell''_0$ and using the second case an isomorphism between fields of $\ell''_0$ and $\ell'_0$. Then the required isomorphism is their composition.

Clearly, if we choose $\ell_0, m_0$ to be the axes of $PP^*_{wu}(\mathcal{F})$ and $Z_0 = (1, 1)$,
then the constructed field will be isomorphic to $\mathcal{F}$.
Hence, it will be isomorphic to $RF^*_{field}(PP^*_{wu}(\mathcal{F}))$ for every choice of $\ell_0, m_0, Z_0$.

This result can also be found in ~\cite[Theorem 5.9]{hall1943projective}.

(iv) Let $\Pi$ be a metric Wu plane and let $\mathcal{F}=RF^*_{field}(\Pi)$. In order to establish an isomorphism between $PP^*_{wu}(\mathcal{F})$ and $\Pi$ we need to define two maps: a map of points and a map of lines.

Points of $PP^*_{wu}(\mathcal{F})$ are pairs $(x,y)\in \mathcal{F}^2$. We recall that the universe of $\mathcal{F}$ is the set of points incident to the line $\ell_0$, and that there is a definable bijection $f$ between the points of $\ell_0$ and $m_0$ (the coordinate axes in $\Pi$, parameters of the considered interpretation). Hence, given $(x,y)$ we can define two auxiliary lines: $m_x$ going through $x\in\ell_0$ and parallel to $m_0$, and $\ell_y$, going through $f(y)\in m_0$ and parallel to $\ell_0$. Let
$A$ be the intersection of $m_x$ and $\ell_y$. We map $(x,y)$ to $A$, and it is clear that $A$ has coordinates $(x,y)$. Thus, we have described a (definable) bijection between the sets of points of $PP^*_{wu}(\mathcal{F})$ and $\Pi$.

Lines of $PP^*_{wu}(\mathcal{F})$ can be specified by equations $ax+by+c=0$ where not all of $a,b,c$ are $0$. Thus, a line is interpreted by a triple $(a,b,c)\in \mathcal{F}^3$. Two lines are defined to be equal if $(a,b,c)$ and $(a',b',c')$ are proportional. A point $(x,y)$ is incident to a line $(a,b,c)$ if $ax+by+c=0$.

Each line in $PP^*_{wu}(\mathcal{F})$ is equal to a line defined by the equation $y=ax+b$ or to a vertical line $x=c$. We construct the corresponding line in $\Pi$ by drawing, in the first case, a line through the point $(0,b)$ with the slope $a$, and in the second case a vertical line (parallel to $m_0$) through $(c,0)$. This maps lines in $PP^*_{wu}(\mathcal{F})$ to lines in  $\Pi$ and is clearly a (definable) bijection preserving the incidence relation.

Concerning the orthogonality relation, we may assume that the coordinate axes $\ell_0$ and $m_0$ in $\Pi$ are selected to be orthogonal. (By Wu, the field $\mathcal{F}$ does not depend on the choice of parameters, up to isomorphism.) We can define
lines $(a,b,c)$ and $(a',b',c')$ in  $PP^*_{wu}(\mathcal{F})$ to be orthogonal iff $aa'+bb'=0$. Then the usual arguments show that this agrees with the orthogonality of the corresponding lines in $\Pi$.

\item[(v)] In the proof of (iii) we use the isomorphisms established by Wu in \cite[page 42, Theorem 1]{bk:Wu1994}. The construction of these isomorphisms mostly involves drawing various parallel lines and is clearly definable. The isomorphism between $\mathcal{F}$ and $RF^*_{field}(PP^*_{wu}(\mathcal{F}))$ where $\ell_0, m_0$ are the axes of $PP^*_{wu}(\mathcal{F})$ and $Z_0 = (1, 1)$, is defined by a map $x \mapsto (x, 0)$, which is first-order definable. The construction given in the proof of (iv) explicitly provides us with a definable isomorphism.

\end{proof}


\section{Undecidability} \label{se:undec}
To establish the undecidability of the theory of Pythagorean fields we refer to the following
theorem of M. Ziegler~\cite{ar:ziegler,BeesonZiegler}:

\begin{theorem} \label{Ziegler}
Let $T$ be a finite subtheory of the theory of the field of reals $(\mathbb{R}; +, \times)$. Then
\begin{enumerate}[(i)]
\item $T$ is undecidable;
\item The same holds for the extension of $T$ by the axioms stating that the characteristic of the field is $0$.
\end{enumerate}
\end{theorem}

Although the second part is not mentioned as a result in ~\cite{ar:ziegler,BeesonZiegler},
it easily follows from Ziegler's proof.

\begin{corollary} \label{Euc}
The theories of Pythagorean fields and of Pythagorean fields of characteristic $0$ are  undecidable.
\end{corollary}

Using the bi-interpretability of Pythagorean fields and metric Wu planes we obtain the undecidability of the theory of metric Wu planes. In fact, we prove a more general theorem establishing the undecidability of a sufficiently wide class of geometric theories. As a preparation for its proof, we define syntactic translations between formulas in the language of fields and formulas in the language $\tau_{wu}$.

Consider any formula $\phi$ in the language of fields. Using the formulas $add$ and $mult$ from the construction of $RF^*_{field}(\Pi)$
to interpret addition and multiplication, we obtain a formula $\phi^{wu}(\ell_0,m_0,Z_0)$ in the language $\tau_{wu}$,
where $\ell_0, m_0, Z_0$ are the parameters (free variables) of the formula.
\newcommand{\Par}{\textit{Par}}

Let $\Par(\ell_0,m_0,Z_0)$ denote the formula stating that $\ell_0$ and $m_0$ are lines
intersecting in exactly one point and that $Z_0$ is not incident with either $\ell_0$ or $m_0$.
These conditions definably specify the admissible values of the parameters. Then, for any metric Wu plane $\Pi$,
\begin{equation}\label{wu}RF^*_{field}(\Pi) \models \phi \Longleftrightarrow \Pi \models \forall{\ell_0,m_0,Z_0}\:(\Par(\ell_0,m_0,Z_0)\to \phi^{wu}).\end{equation}

Similarly, for the other interpretation, consider any formula $\phi$ in the language of $\tau_{wu}$. Using the formulas from the construction of $PP^*_{wu}$ to interpret the two sorts of variables, equality, incidence and orthogonality, we obtain a formula $\phi^{field}$ in the language of fields. Then, for any field $\mathcal{F}$,
\begin{equation}\label{field}PP^*_{wu}(\mathcal{F}) \models \phi \Longleftrightarrow \mathcal{F} \models \phi^{field}.\end{equation}

Now we are ready to state a geometric version of Ziegler's Theorem.

Let $\Pi_\mathbb{R}$ denote the real plane $PP^*_{wu}(\mathbb{R})$. Let $\WU$ denote the first order theory of metric Wu planes and let  $\PF$ denote the first order theory of Pythagorean fields of characteristic $0$.

\begin{theorem} \label{G-Ziegler}
Let $T$ be  a finite set of axioms in the vocabulary $\tau_{wu}$ such that $\Pi_\mathbb{R} \models T$. Then $T \cup \WU$  is undecidable.
\end{theorem}

\begin{proof}
Let $T' = \{ \phi^{field} \mid \phi \in T\}$. Then by Ziegler's theorem, $T' \cup \PF$ is undecidable.

We want to prove that
\begin{equation}\label{eq} T' \cup \PF \models \phi \iff  T \cup \WU \models \forall{\ell_0,m_0,Z_0}\:(\Par(\ell_0,m_0,Z_0)\to \phi^{wu}).
\end{equation}
Then, since the translation $(\cdot)^{wu}$ is computable,
this provides a computable reduction of $T' \cup \PF$ to $T \cup \WU$ and proves that the latter is undecidable.

To prove (\ref{eq}), suppose $T \cup \WU \models \forall{\ell_0,m_0,Z_0}\:(\Par(\ell_0,m_0,Z_0)\to \phi^{wu})$. Take any $\mathcal{F} \models T' \cup \PF$. Then by Theorem \ref{Mak}, $PP^*_{wu}(\mathcal{F})$ is a metric Wu plane and, using (\ref{field}), we have $PP^*_{wu}(\mathcal{F}) \models T \cup \WU$. Hence, $$PP^*_{wu}(\mathcal{F}) \models \forall{\ell_0,m_0,Z_0}\:(\Par(\ell_0,m_0,Z_0)\to \phi^{wu}).$$ Then, by (\ref{wu}), $RF^*_{field}(PP^*_{wu}(\mathcal{F})) \models \phi$ and, since ${\mathcal{F} \cong RF^*_{field}(PP^*_{wu}(\mathcal{F}))}$, we obtain $T' \cup \PF \models \phi$.

Suppose $T' \cup \PF \models \phi$. Consider  any $\Pi \models T \cup \WU$. Let $\mathcal{F} = RF^*_{field}(\Pi)$. By Theorem \ref{Mak}, $\mathcal{F}$ is a Pythagorean field of characteristic $0$ and $\Pi$ is isomorphic to $PP^*_{wu}(\mathcal{F})$. Then, using (\ref{field}), we obtain  $\mathcal{F} \models T' \cup \PF$. Therefore, $\mathcal{F} \models \phi$ and by (\ref{wu}), $\Pi \models \forall{\ell_0,m_0,Z_0}\:(\Par(\ell_0,m_0,Z_0)\to \phi^{wu}).$ It follows that $T \cup \WU \models \forall{\ell_0,m_0,Z_0}\:(\Par(\ell_0,m_0,Z_0)\to \phi^{wu})$.

This completes the proof of (\ref{eq}) and thereby of Theorem~\ref{G-Ziegler}.
\end{proof}

\begin{corollary} \label{undecidability}
The theory of metric Wu planes is undecidable.
\end{corollary}

%% file: hjaxioms-rev.tex
\section{Logical Huzita -- Justin axioms}
\label{se:hjaxioms}

Looking at the Huzita -- Justin axioms it appears as if one could axiomatize
origami geometry in the language of incidence and orthogonality only,
taking as the underlying geometry a metric Wu plane.
However, it turns out that there are multiple reasons for choosing
as the underlying geometry an {\em ordered} metric Wu plane.

Firstly, the Huzita -- Justin axiom (H-3), which says that for any two lines $l_1, l_2$ there is a fold which places
$l_1$ on $l_2$, is only true if $l_1, l_2$ are  non-isotropic (not orthogonal to themselves).
Our formulation of (H*3) takes this into account.

Secondly, Axiom (H-5) states that, given two points $P_1$ and $P_2$ and a line $\ell$, one should be able to construct a fold that places $P_1$ onto $\ell$ and passes through $P_2$. However, such a fold only exists provided $P_1$ is closer to the given line $\ell$ than $P_2$. This condition needs to be expressible in the language, which we achieve by introducing the betweenness relation, as formulated in (H*5). Other possibilities are discussed at the end of the paper.

We say that a metric Wu plane is orderable if it can be equipped with a ternary relation $Be(A,B,C)$ which satisfies the Hilbertian axioms of betweenness.
In Proposition \ref{prop:orderable-2} 
we prove that a metric Wu plane is  orderable iff it has no isotropic lines.

If $A,B,C$ are collinear on a line $\ell$, let  $\ell'$ be orthogonal to $\ell$ going through the point $C$ and
let $B'$ be the point on $\ell$ obtained by placing $B$ on $\ell$ after folding along $\ell'$.
In the real plane and the presence of the betweenness relation $B$ is between $A$ and $C$, $Be(A,B,C)$
we have $Out(A,C,B)$ iff either $Be(A,B,C)$ or $Be(A,B',C)$.
Therefore $Out(A,C,B)$ is definable using $Be(A,B,C)$ and the usual axioms for betweenness.
It is not obvious, however, if $Be(A,B,C)$ can be formulated using an axiomatization of $Out(A,C,B)$.
The same can be said about $Closer(A_1, A_2, \ell)$.

A reasonable axiomatization of origami geometry can be obtained from ordered metric Wu planes by adding a finite set of axioms.

Metric Wu planes already satisfy Huzita -- Justin axioms (H-1), (H-2), (H*3), (H-4)
where (H*3) is our modification of (H-3), see Proposition \ref{prop:isotropic}.
An {\em ordered origami plane} is an ordered metric Wu plane which satisfies also our modified axioms (H*5) and (H*6), see section \ref{se:vieta}.
This deviates from the definition given in \cite{makowsky2018undecidability,makowsky2019can}
where an ordered origami plane is defined in an inconsistent way.

The ordered metric Wu planes are bi-interpretable with ordered Pythagorean fields, whose first-order theory is undecidable by Ziegler's theorem. Similarly, we obtain that the first order theory of our axiomatization of origami geometry is also undecidable.

Huzita--Justin axioms were not meant to be axioms in the logical sense, but rather rules of folding.
Yet, one can try to naively formulate them in the language $\tau_{wu}$ by treating the
requirement to construct an object (satisfying given conditions) by a classical existential statement.
Huzita--Justin axioms are naturally stated using the relation $\Sym(P_1, \ell, P_2)$
\emph{``points $P_1$ and $P_2$ are symmetric with respect to line $\ell$''} defined
in Section~\ref{se:wuplanes}. Then one obtains the following versions of Huzita--Justin axioms.

\begin{description}
\item[(H-1):]
Given two points $P_1$ and $P_2$, one can make a fold that passes through both of them:
$$
\forall P_1, P_2 \,\exists \ell \:(P_1 \in \ell \wedge P_2 \in \ell).
$$
\item[(H-2):]
Given two points $P_1$ and $P_2$, one can make a fold that places $P_1$ onto $P_2$:
$$
\forall P_1, P_2\, \exists \ell\: \Sym(P_1, \ell, P_2).
$$
\item[(H-3):]
Given two lines $\ell_1$ and $\ell_2$, one can make a fold that places $\ell_1$ onto $\ell_2$:
$$
\forall \ell_1, \ell_2\, \exists k\, \forall P_1\:
\left(
P_1 \in \ell_1 \rightarrow \exists P_2
\left( P_2 \in \ell_2 \wedge \Sym(P_1, k, P_2) \right)
\right).
$$

\item[(H-4):]
Given a point $P$ and a line $\ell$, one can make a fold orthogonal to $\ell$ that passes through $P$:
$$
\forall P, \ell\, \exists k\: (P \in k \wedge \ell \perp k).
$$
\item[(H-5):]
Given two points $P_1$ and $P_2$ and a line $\ell_1$, one can make a fold that places $P_1$ onto $\ell_1$ and passes through $P_2$:
$$
\forall P_1, P_2, \ell_1 \,\exists \ell_2\: (P_2 \in \ell_2 \wedge \exists P_3\: (\Sym(P_1, \ell_2, P_3) \wedge P_3 \in \ell_1)).
$$
\item[(H-6):]
Given two points $P_1$ and $P_2$ and two lines $\ell_1$ and $\ell_2$, one can make a fold that places $P_1$ onto $\ell_1$ and $P_2$ onto $\ell_2$:
\begin{multline*}
\forall P_1, P_2, \ell_1, \ell_2 \, \exists \ell_3 \:
\left(
\exists Q_1\: (\Sym(P_1, \ell_3, Q_1) \wedge Q_1 \in \ell_1) \wedge
\right.
\\
\left.
\wedge
\exists Q_2\: (\Sym(P_2, \ell_3, Q_2) \wedge Q_2 \in \ell_2)
\right).
\end{multline*}
\end{description}

Since the original Huzita -- Justin axioms talk only about the possibility of the existence of folds,
Axioms (H-5) and (H-6) formulated above do not hold in a real plane. The exceptional configurations, where a described fold does
not exist, have to be described explicitly.
In order to fix that we are going to amend them with the most obvious conditions under which the lines would indeed exist.
On the other hand, the formalizations of Axioms (H-1, H-2, H-3, H-4) obviously hold in the real plane and, as explained below, almost hold in any metric Wu plane.

We would like to state that metric Wu planes satisfy the first four origami axioms (H-1, H-2, H-3, H-4). However, since there could exist non-isotropic lines in a metric Wu plane, they do not necessarily satisfy (H-3). Thus, we modify this axiom.

\begin{description}
\item[(H*3):]
Given two non-isotropic lines $\ell_1$ and $\ell_2$, there is a fold (line) that places $\ell_1$ onto $\ell_2$.
\end{description}

\begin{proposition}
\label{prop:isotropic}
Every metric Wu plane satisfies the origami axioms (H-1, H-2, H*3, H-4).
\end{proposition}

\begin{proof}
(H-1) is equivalent to (I-1) and (H-4) is equivalent to (O-2). To prove (H-2) we use the construction from \cite[page 75]{bk:Wu1994}.

(AxSymAx) is an analogue of (H*3) for intersecting lines, so we may only consider the case of parallel lines. Take $\ell_1$ and $\ell_2$ parallel to each other. Take any point $P_1 \in \ell_1$ and drop a perpendicular from $P_1$ on $\ell_2$. Let the intersecting point be $P_2$. Then we claim that the perpendicular bisector of $P_1 P_2$ is the line we need.
\end{proof}

%% file: ordered-rev.tex
\section{Ordered metric Wu planes and Pythagorean fields}
\label{se:ordered}

We have established a correspondence between Pythagorean fields and metric Wu planes. Following \cite{alperin2000mathematical}, the next step would be a correspondence between Euclidean fields, defined below, and planes satisfying some analogue of (H-5).

\begin{definition}
\label{eucfield}
A \emph{Euclidean field} is a formally real Pythagorean field such that every element is either a square or the  opposite of a square:
$$\forall x \exists y \ (x=y^2 \lor -x=y^2).$$
\end{definition}

As we will discuss in the next section, Euclidean fields are always uniquely ordered, therefore we want our plane to be in some sense ``ordered'' as well. One way to do so would be to take as an axiom that there are no isotropic lines. Then the corresponding field would be formally real and hence orderable (see below).

We take a different approach and follow \cite{bk:Wu1994}. As the concept 'lying between' is not definable in a metric Wu plane, we introduce a new relation of Betweenness and the axioms that describe it. We interpret $Be(P_1, P_2, P_3)$ as \emph{three distinct points are on the same line and $P_2$ is between $P_1$ and $P_3$}. Note, that if three points are not all distinct, none of them lie between the others. Let $\tau_{o-wu}$ be the signature consisting of $\in$, $\perp$ and $Be$.

\paragraph{Axioms of betweenness} 
\begin{description}
\item[(B-1):]
Let $A, B, C$ be three distinct points on a line. If $B$ lies between $A$ and $C$, then $B$ also lies between $C$ and $A$.
\item[(B-2):]
For any two distinct points $A$ and $C$ on a line, there always exists another point $B$ which lies between $A$ and $C$, and another point $D$ such that $C$ lies between $A$ and $D$. 
\item[(B-3):]
Given any three distinct points $A, B, C$ on a line, one and only one of the following three cases holds: $B$ lies between $A$ and $C$, $A$ lies between $B$ and $C$, and $C$ lies between $A$ and $B$.
\item[(B-4):] (Pasch)
Assume the points $A, B, C$ and $\ell$ in general position,
i.e.
the three points are not on one line, none of the points is on $\ell$.
Suppose there is a point $D$ at which $\ell$ and the line $AB$ intersect.
If $Be(A, D, B)$ there is $D' \in \ell$ with
$Be(A, D', C)$ or
$Be(B, D', C)$.
\end{description}

\begin{definition}
A $\tau_{o-wu}$ structure $\Pi$ is an ordered metric Wu plane if it is a metric Wu plane satisfying axioms of betweenness (B-1, B-2, B-3, B-4).
\end{definition}

\begin{proposition}
\label{prop:orderable-1}
Every ordered metric Wu plane satisfies (H-3) and hence the origami axioms (H-1, H-2, H-3, H-4).
\end{proposition}

\begin{proof}
There are no isotropic lines in ordered metric Wu planes as proven in \cite[page 107, Theorem 3]{bk:Wu1994}.
\end{proof}

If $\mathcal{F}$ is an ordered field, we define the relation of betweenness on $PP^*_{wu}(\mathcal{F})$ in the standard way. If $\Pi$ is an ordered metric Wu plane, we follow \cite[page 105]{
bk:Wu1994} to define an order on $RF^*_{field}(\Pi)$. Using the result \cite[page 103, Separation property 1]{bk:Wu1994}, we separate all points on $\ell_0$ distinct from $0$ into two parts. Then $0$ lies between $A, B$ when $A, B$ lie in different parts, and $0$ does not lie between $A, B$ when $A, B$ lie in the same part. We define those numbers in $RF^*_{field}(\Pi)$ whose corresponding points lie in the same part as $1$ to be positive numbers and those whose corresponding points lie in the other part to be negative numbers. Then we can say that $a < b$ whenever $b - a$ is a positive number.

\begin{theorem} \label{o-wuplanes}
\begin{itemize}
\item[(i)] Let $\mathcal{F}$ be an ordered Pythagorean field. 
Then $PP^*_{wu}(\mathcal{F})$ is an ordered metric Wu plane.
\item[(ii)] Let $\Pi$ be an ordered metric Wu plane. Then $RF^*_{field}(\Pi)$ is an ordered Pythagorean field.
\item[(iii)] $RF^*_{field}(PP^*_{wu}(\mathcal{F}))$ is isomorphic to $\mathcal{F}$.
\item[(iv)] $PP^*_{wu}(RF^*_{field}(\Pi))$ is isomorphic to $\Pi$.

\item[(v)] The isomorphisms in (iii) and (iv) are definable, that is, form a bi-interpretation between the classes of ordered metric Wu planes and of ordered Pythagorean fields. 
\end{itemize}
\end{theorem}

\begin{proof}
(i) Using properties of ordered fields, it is easy to check that $PP^*_{wu}(\mathcal{F})$ satisfies the axioms of betweenness (B-1, B-2, B-3, B-4).

(ii) We only need to show that $RF^*_{field}(\Pi)$ is an ordered field. This is proved on \cite[page 105, Theorem 1]{bk:Wu1994}.

(iii) By Theorem \ref{Mak} it is sufficient to check that the relation of order is preserved. As discussed on \cite[page 105]{bk:Wu1994}, if we take different $\ell_0$, $m_0$, $Z_0$ in the construction of $RF^*_{field}(\Pi)$, the canonical isomorphism between the obtained fields will preserve order. If we choose $\ell_0, m_0$ to be the axes of $PP^*_{wu}(\mathcal{F})$ and $Z_0 = (1, 1)$, then the field is clearly isomorphic to $\mathcal{F}$. It follows that $\mathcal{F}$ is isomorphic to  $RF^*_{field}(PP^*_{wu}(\mathcal{F}))$ for any choice of parameters.

(iv) Let $\mathcal{F}$ denote the field  $RF^*_{field}(\Pi)$. By Theorem \ref{Mak} it is sufficient to check that the betweenness relation is preserved under the isomorphism of metric Wu planes $\Pi$ and $PP^*_{wu}(\mathcal{F})$. 

Since the collinearity is preserved, it is sufficient to consider the betweenness relation for points on the same line. Suppose three points $A,B,C$ on a line $\ell$ in $\Pi$ are given. Assume $\ell$ is not vertical and consider the intersections $A', B', C'$ of $\ell_0$ axis and lines parallel to $m_0$ and going through $A, B, C$ respectively. By Corollary 1 in \cite[page 104]{bk:Wu1994}, $Be(A, B, C)$ holds in an ordered metric Wu plane if and only if $Be(A', B', C')$ does. On the other hand, the interpretation of betweenness in $PP^*_{wu}(\mathcal{F})$ for points of the coordinate axis is the same as that in $\Pi$. This shows the claim in the case $\ell$ is not vertical. 

If $\ell$ is vertical, we consider the projections of $A,B,C$ on the $m_0$ axis and the corresponding points on $\ell_0$ via the bijection $f$. By the same principle, $f$ preserves the betweenness on respective coordinate axes, hence $Be(A,B,C)$ holds in $\Pi$ iff it holds in $PP^*_{wu}(\mathcal{F})$.

\item[(v)] We use the same isomorphisms as in Theorem \ref{Mak}, which are already showed to be definable.
\end{proof}

After establishing Theorem \ref{o-wuplanes}, we also want to obtain an ordered analogue of Theorem \ref{G-Ziegler}. For this purpose, we need an ordered version of Theorem \ref{Ziegler}. Although in ~\cite{ar:ziegler,BeesonZiegler} only theories in the language of fields are concerned, the proof of Ziegler's Theorem essentially consists of constructing a model of $T$, in which the ring of integers is interpretable. Then it is easy to see that the proof still holds for the case of ordered fields, which gives us the following result.

\begin{theorem} \label{o-Ziegler}
Let $T$ be a finite subtheory of the theory of the ordered field of reals $(\mathbb{R}; +, \times, \leq)$. Then 
\begin{enumerate}[(i)] 
\item $T$ is undecidable; 
\item The same holds for the extension of $T$ by the axioms stating that the characteristic of the field is $0$.  
\end{enumerate}
\end{theorem}

Then using Theorem \ref{o-wuplanes} and the same technique as in Theorem \ref{G-Ziegler}, we obtain an ordered version of the geometrical Ziegler's theorem.

Let $\Pi_\mathbb{R}$ denote the real plane $PP^*_{wu}(\mathbb{R})$. Let $\OWU$ denote the first order theory of ordered metric Wu planes.

\begin{theorem} \label{o-G-Ziegler}
Let $T$ be a finite set of axioms in the vocabulary $\tau_{o-wu}$ such that $\Pi_\mathbb{R} \models T$. Then $T \cup \OWU$  is undecidable.
\end{theorem}


Next we consider the question of orderability. 
Recall that a field is orderable (or \emph{formally real}) if $-1$ is not a sum of squares. For Pythagorean fields this is equivalent to saying that $-1$ is not a square. The statement that there are no isotropic lines plays a similar role for metric Wu planes. 

\begin{proposition}
\label{prop:orderable-2}
A metric Wu plane $\Pi$ is orderable iff there are no isotropic lines in $\Pi$.
\end{proposition}

\begin{proof}
In one direction, we have already mentioned a theorem of Wu that ordered metric Wu planes have no isotropic lines. In the other direction, we assume a $\Pi$ without isotropic lines is given and consider as a coordinate system a pair of orthogonal lines and the corresponding field  ${\cal F}=RF^*_{field}(\Pi)$. By Theorem \ref{Mak} (iv)  $PP^*_{wu}({\cal F})$ is isomorphic to $\Pi$. 

We claim that ${\cal F}$ is formally real. Assume the contrary, then $d^2+1=0$ in $\cal F$, for some $d$. Then the line defined by the points $(1,0)$ and $(0,d)$ (and the parallel line given by the equation $dx+y=0$) is isotropic. 
Since $\cal F$ is formally real, $PP^*_{wu}({\cal F})$ is orderable, but it is isomorphic to $\Pi$.
\end{proof}


%% file: vieta-rev.tex
\section{Euclidean and Vieta fields and origami axioms}
\label{se:vieta}

As already mentioned in the previous section, our goal is to establish a correspondence between Euclidean fields and planes satisfying some amended version of (H-5). Recall that a Euclidean field is a formally real Pythagorean field such that every element is either a square or the  opposite of a square.

The nonzero squares of a Euclidean field constitute a positive cone, hence (see~\cite{Becker}) Euclidean fields admit a unique ordering: 
$$x\leq y \iff \exists z\ (x+z^2=y).$$ Since the ordering is definable, one often  considers Euclidean fields as ordered fields. An ordered field is Euclidean iff each positive element in it is a square. 

\begin{proposition}
The first order theory of Euclidean fields is undecidable.
\end{proposition}

\begin{proof}
Ziegler's Theorem.
\end{proof}

Next we formulate our amended version of Axiom (H-5) in which we add an appropriate precondition for the 
constructed fold to exist. 
Below we use the notion \emph{a point $A$ is closer to a line $\ell$ than to a point $B$}, $\textit{Closer}(A, \ell, B)$,
which can be formulated in the language $\tau_{o-wu}$ by saying that 
there exist points $H \in \ell$ and $B'$ and a line $m\ni A$ such that $\Sym(H,m,B')$, while $Be(A, B', B)$ or $B = B'$.

\begin{description}
\item[(H*5):]
Given two points $P_1$ and $P_2$ and a line $\ell_1$, if $P_2$ is closer to $\ell_1$ than to $P_1$, 
then there is a fold (line) that places $P_1$ onto $\ell_1$
and passes through $P_2$.
\end{description}



\begin{definition}
A $\tau_{o-wu}$ structure $\Pi$ is a \emph{Euclidean ordered metric Wu plane} if it is an ordered metric Wu plane satisfying (H*5).
\end{definition}

\begin{theorem} \label{th:pyth}
\begin{itemize}
\item[(i)] Let $\mathcal{F}$ be a Euclidean field. 
Then $PP^*_{wu}(\mathcal{F})$ is a Euclidean ordered metric Wu plane.
\item[(ii)] Let $\Pi$ be a Euclidean ordered metric Wu plane.
Then $RF^*_{field}(\Pi)$ is a Euclidean field.
\item[(iii)] Furthermore, $RF^*_{field}(PP^*_{wu}(\mathcal{F}))$ is isomorphic to $\mathcal{F}$.
\item[(iv)] $PP^*_{wu}(RF^*_{field}(\Pi))$ is isomorphic to $\Pi$.
\item[(v)] The isomorphisms in (iii) and (iv) are definable, that is, form a bi-interpretation between the classes of Euclidean metric Wu planes and of Euclidean fields.
\end{itemize}
\end{theorem}

\begin{proof}
\begin{itemize}
    \item[(i)]  In order to show (H*5) it is enough to prove that a circle intersects a line whenever the radius is smaller than the distance between the center and the line. In $PP^*_{wu}(\mathcal{F})$ that is equivalent to solving a quadratic equation. Therefore, if $\mathcal{F}$ is Euclidean, (H*5) holds.

    \item[(ii)] Consider any positive $s \in RF^*_{field}(\Pi)$ and show that the square root of $s$ exists. Without loss of generality we assume $s > 1$, since otherwise we can find a square root of $s^{-1}$. Let $P_1 = (0, s)$, $P_2 = (0, \frac{s - 1}{2})$ and let $\ell_0$ be the $x$ axis. Then by (H*5) we can find a point $P_3 = (x, y)$, such that $P_3 \in \ell_0$ and $|P_1 P_2| = |P_3 P_2|$. The first condition gives us $y = 0$. The second one, if we calculate the squares of distances, means ${\left( \frac{s-1}{2} \right)^2 + x^2 = \left( \frac{s+1}{2} \right)^2}$, which is equivalent to $x^2 = s$.
    
    \item[(iii-v)] We use Theorem \ref{o-wuplanes}, since every Euclidean field is a Pythagorean field and every Euclidean ordered metric Wu plane is an ordered metric Wu plane.
\end{itemize}
\end{proof}


Then by Theorem \ref{o-G-Ziegler}, we obtain:

\begin{corollary}
The theory of Euclidean ordered Wu planes is undecidable.
\end{corollary}

Finally, we would like to find a correct version of (H*6) and to establish its correspondence with Vieta fields as defined below.

\begin{definition}
A \emph{Vieta field} is a Euclidean field in which every element is a cube:
$$\forall x \ \exists y \ y^3 = x.$$
\end{definition}

It follows from Cardano formula that any cubic polynomial over a Vieta field has at least one root.
\begin{proposition}
The first order theory of Vieta fields is undecidable.
\end{proposition}

\begin{proof}
Ziegler's Theorem.
\end{proof}

The following version of (H-6) is inspired by~\cite[Proposition 6]{ghourabi2012algebraic}.

\begin{description}
\item[(H*6):]
Given two points $P_1$ and $P_2$ and two lines $\ell_1$ and $\ell_2$, if $P_1 \notin \ell_1$, $P_2 \notin \ell_2$, $\ell_1$ and $\ell_2$ are not parallel and points are distinct or lines are distinct, then there is a fold (line) that places $P_1$ onto $\ell_1$ and $P_2$ onto $\ell_2$.\\ 
\end{description}



\begin{definition}
$\Pi$ is an \emph{ordered Wu origami plane} if it is an ordered metric Wu plane which also satisfies (H*5) and (H*6).
\end{definition}

\begin{theorem} \label{th:vieta}

\begin{itemize}
\item[(i)] Let $\mathcal{F}$ be a Vieta field. 
Then $PP^*_{wu}(\mathcal{F})$ is an ordered Wu origami plane.
\item[(ii)] Let $\Pi$ be an ordered Wu origami plane. Then $RF^*_{field}(\Pi)$ is a Vieta field.
\item[(iii)] Furthermore, $RF^*_{field}(PP^*_{wu}(\mathcal{F}))$ is isomorphic to $\mathcal{F}$.
\item[(iv)] $PP^*_{wu}(RF^*_{field}(\Pi))$ is isomorphic to $\Pi$.
\item[(v)] The isomorphisms in (iii) and (iv) are definable, that is, form a bi-interpretation between the classes of ordered Wu origami planes and of Vieta fields.
\end{itemize}
\end{theorem}

\begin{proof}

\begin{itemize}
\item[(i)]
It sufficies to show that (H*6) holds in $PP^*_{wu}(\mathcal{F})$. We use ~\cite[Proposition 6]{ghourabi2012algebraic} to conclude that the conditions we chose are sufficient for the existence of a fold. Although the original result was proven specifically for the real plane, we note that it essentially uses only the Vieta property of $\mathbb{R}$ and therefore holds for any plane over a Vieta field.

\item[(ii)]  Take any $r \in RF^*_{field}(\Pi)$. Let $P_1 = (-1, 0)$, $P_2 = (0, -r)$, $\ell_1$ be the line $x = 1$ and $\ell_2$ the line $y = r$. Then by (H*6) we can find a line $\ell$. Drop perpendiculars from $P_1$ and $P_2$ on $\ell$ and let the constructed points be $H_1$ and $H_2$. Then since a reflection in $\ell$ maps $P_1$ onto $\ell_1$ and $P_2$ onto $\ell_2$, the coordinates of $H_1$ and $H_2$ have to be of the form $(0, s)$ and $(t, 0)$ respectively. Lines $P_1 H_1$ and $P_2 H_2$ have to be parallel, therefore $t = \frac{r}{s}$. Finally, $H_1 H_2$ is perpendicular to $P_1 H_1$, which means $\frac{r}{s} - s^2 = 0$. Hence, $s^3 = r$.
\item[(iii-v)] Once again these follow from Theorem \ref{o-wuplanes}, since every Vieta field is a Pythagorean field and every ordered Wu origami plane is an ordered metric Wu plane.
\end{itemize}
\end{proof}

\begin{corollary}
The theory of ordered Wu origami planes is undecidable.
\end{corollary}

%% file: jam-dconclu.tex
\section{Discussion}
\label{se:conclu}


We have axiomatized the classes of orthogonal Wu planes using versions of origami axioms and established their bi-interpretations with the first order theories of fields (corresponding to the classes of Pythagorean, Euclidean and Vieta fields). A few natural questions concerning the axiomatization of geometry via origami constructions were left open. 

One such question concerns the choice of the considered language. Although the orthogonality of lines is natural from the point of view of origami constructions --- orthogonality can be tested just by a single fold --- we see that the Huzita--Justin axioms are easily formulated using the notion of symmetry of two points w.r.t.\ a line $\Sym(P,\ell,Q)$. It would be natural to  
consider this notion as basic and orthogonality as definable. The $\Sym$ predicate behaves well provided the metric Wu plane is orderable, that is, has no isotropic lines. 

\begin{problem}
Find a natural axiomatization of orderable  metric Wu planes in terms of $\in$ and $\Sym$. 
\end{problem}

A similar question can be asked about betweenness. We have based our axiomatization on the standard Hilbertian axioms for betweenness. One can, however, consider as basic the relation $\textit{Closer}(A,\ell,B)$ which holds if $A$ is closer to line $\ell$ than to $B$. It is the relation that was used in the statement of (H*5). 

\begin{problem} Find a natural axiomatization of the class of ordered metric Wu planes in terms of $\in$, $\perp$ and \textit{Closer}. In particular, this requires that there is a first-order formula that works as a definition of betweenness in each structure satisfying these axioms.
\end{problem}

Another question concerns the definability of betweenness in Euclidean metric Wu planes. 
Recall that in Euclidean fields the ordering is definable. This suggests that there is an axiomatization of Euclidean ordered metric Wu planes in the language $\tau_{wu}$ only. In fact, one such axiomatization based on the so-called \emph{Euclidean axiom of betweenness} is well-known~\cite[page 149]{degen-profke}. This axiom defines the betweenness relation $Be(A,B,C)$ by stating that $A,B,C$ are collinear and there exists a point $D$ such that $DA\perp DC$ and $DB\perp AC$. Then, a metric Wu plane is Euclidean ordered iff it has no isotropic lines and the above relation $Be$ satisfies the usual axioms of betweenness. 

It would be interesting to know if the Euclidean axiom of betweenness can be replaced by its alternative suggested by the origami axiom (H-5). First, \emph{define} $\textit{Closer}(A,\ell,B)$ by saying that there is a fold $m$ that goes through $A$ and places $B$ on $\ell$. Second, define $Be(A,B,C)$ by saying that  $A,B,C$ are collinear and there is a line $\ell\ni D$ such that $\ell\perp AC$, $\textit{Closer}(A,\ell,C)$ and $\textit{Closer}(C,\ell,A)$. State (some of) the betweenness axioms for the defined relation. Though this approach may work, there are a number of details to be worked out here that we leave for a future study. 

\begin{problem}
Find an axiomatization of Euclidean orderable Wu planes in the language $\tau_{wu}$ that would be natural from the point of view of origami.
\end{problem}


Yet another interesting direction of study is to develop a constructive version of origami geometry as a logical theory based on intuitionistic logic, in the spirit of the work of Beeson~\cite{Beeson}. In such a theory existential statements would yield actual origami constructions rather than just be classically true.  

\begin{problem}
Develop a constructive version of origami geometry. 
\end{problem}

\section{Acknowledgements} 

The work of Lev Beklemishev and Anna Dmitrieva was supported by the Academic Fund Program at the National Research University Higher School of Economics (HSE) in 2019--2020 (grant No. 19-04-050) and by the Russian Academic Excellence Project ``5--100''.

%% file: jam-rev.bbl
\begin{thebibliography}{GKK12}

\bibitem[Alp00]{alperin2000mathematical}
R.C. Alperin.
\newblock A mathematical theory of origami constructions and numbers.
\newblock {\em New York J. Math}, 6(119):133, 2000.

\bibitem[Bac73]{bk:Bach1973}
F.~Bachmann.
\newblock {\em Aufbau der Geometrie aus dem Spiegelungsbegriff}.
\newblock Springer-Verlag, 1973.

\bibitem[Bec74]{Becker}
Eberhard Becker.
\newblock Euklidische {Körper und euklidische Hüllen von Körpern}.
\newblock {\em Journal für die reine und angewandte Mathematik},
  1974(268--269):41--52, 1974.

\bibitem[Bee]{BeesonZiegler}
M.~Beeson.
\newblock Some undecidable field theories.
\newblock Translation of \cite{ar:ziegler}. Available at
  www.michaelbeeson.com/research/papers/Ziegler.pdf.

\bibitem[Bee15]{Beeson}
Michael Beeson.
\newblock A constructive version of {T}arski's geometry.
\newblock {\em Annals of Pure and Applied Logic}, 166(11):1199 -- 1273, 2015.

\bibitem[Bel34]{Beloch}
M.P. Beloch.
\newblock Alcune applicazioni del metodo del ripiegamento della carta di
  sundara-row.
\newblock {\em Atti dell’Acc. di Scienze, Mediche. Nat Mat Ferrara Serie II},
  11:186--189, 1934.

\bibitem[Blu80]{blumenthal1980modern}
L.~M. Blumenthal.
\newblock {\em A modern view of geometry}.
\newblock Courier Corporation, 1980.

\bibitem[DO07]{bk:DemRou}
Erik~D. Demaine and Joseph O'Rourke.
\newblock {\em Geometric folding algorithms}.
\newblock Cambridge University Press, 2007.

\bibitem[DP76]{degen-profke}
W.~Degen and L.~Profke.
\newblock {\em Grundlagen der affinen und euklidischen Geometrie}.
\newblock Teubner, 1976.

\bibitem[Fri18]{Fri-Orig}
M.~Friedman.
\newblock {\em A History of Folding in Mathematics: Mathematizing the Margins}.
\newblock Birkh\"{a}user, 2018.

\bibitem[FV14]{FrVi}
H.~Friedman and A.~Visser.
\newblock When bi-interpretability implies synonymy.
\newblock Logic Group Preprint Series 320, 2014.

\bibitem[GKK12]{ghourabi2012algebraic}
Fadoua Ghourabi, Asem Kasem, and Cezary Kaliszyk.
\newblock Algebraic analysis of {Huzita’s} origami operations and their
  extensions.
\newblock In {\em International Workshop on Automated Deduction in Geometry},
  pages 143--160. Springer, 2012.

\bibitem[Hal43]{hall1943projective}
M.~Hall.
\newblock Projective planes.
\newblock {\em Transactions of the American Mathematical Society},
  54(2):229--277, 1943.

\bibitem[Har00]{bk:Hartshorne2000}
R.~Hartshorne.
\newblock {\em Geometry: Euclid and Beyond}.
\newblock Springer, 2000.

\bibitem[Hil71]{hilbert1971foundations}
D.~Hilbert.
\newblock {\em Foundations of Geometry, Second Edition, translated from the
  Tenth Edition, revised and enlarged by Dr Paul Bernays}.
\newblock The Open Court Publishing Company, La Salle, Illinois, 1971.

\bibitem[Hil13]{hilbert2013grundlagen}
D.~Hilbert.
\newblock {\em Grundlagen der {G}eometrie}.
\newblock Springer-Verlag, 2013.

\bibitem[Hod93]{bk:Hodges93}
W.~Hodges.
\newblock {\em Model Theory}, volume~42 of {\em Encyclopedia of Mathematics and
  its Applications}.
\newblock Cambridge University Press, 1993.

\bibitem[Huz89a]{Huzita}
H.~Huzita.
\newblock Axiomatic development of origami geometry.
\newblock In {\em \emph{\cite{Proc-Huzita}}}, pages 143--158, 1989.

\bibitem[Huz89b]{Proc-Huzita}
H.~Huzita, editor.
\newblock {\em Proceedings of the 1st international meeting of origami, science
  and technology}.
\newblock Comune di Ferrara and Centro Origami Diffusion, Ferrara, 1989.

\bibitem[Iva16]{ivanov2016affine}
N.V. Ivanov.
\newblock Affine planes, ternary rings, and examples of non-desarguesian
  planes.
\newblock {\em arXiv preprint arXiv:1604.04945}, 2016.

\bibitem[Jus89]{justin1989resolution}
J.~Justin.
\newblock R{\'e}solution par le pliage de {\'e}quation du troisieme degr{\'e}
  et applications g{\'e}om{\'e}triques.
\newblock In {\em Proceedings of the first international meeting of origami
  science and technology}, pages 251--261. Ferrara, Italy, 1989.

\bibitem[LA09]{AlpLang}
Robert~J. Lang and Roger~C. Alperin.
\newblock One-, two-, and multi-fold origami axioms.
\newblock In {\em Origami4: Fourth International Meeting of Origami Science,
  Mathematics, and Education}, pages 383--406, 2009.

\bibitem[Lan03]{LangArt}
Robert~J. Lang.
\newblock {\em Origami Design Secrets: Mathematical Methods for an Ancient
  Art}.
\newblock A K Peters, 2003.

\bibitem[Mak04]{ar:MakowskyTARSKI}
J.A. Makowsky.
\newblock Algorithmic uses of the {F}eferman--{V}aught theorem.
\newblock {\em Annals of Pure and Applied Logic}, 126.1-3:159--213, 2004.

\bibitem[Mak17]{makowsky2017can}
Johann~A. Makowsky.
\newblock Can one design a geometry engine? on the (un) decidability of affine
  {E}uclidean geometries.
\newblock {\em arXiv preprint arXiv:1712.07474}, 2017.

\bibitem[Mak18]{makowsky2018undecidability}
Johann~A. Makowsky.
\newblock The undecidability of orthogonal and origami geometries.
\newblock In {\em International Workshop on Logic, Language, Information, and
  Computation}, pages 250--270. Springer, 2018.

\bibitem[Mak19]{makowsky2019can}
Johann~A. Makowsky.
\newblock Can one design a geometry engine?
\newblock {\em Annals of Mathematics and Artificial Intelligence},
  85(2-4):259--291, 2019.

\bibitem[Rao17]{Rao}
T.~Sundara Rao.
\newblock {\em Geometric Exercises in Paper Folding}.
\newblock The Open Court Publishing Company, 1917.
\newblock Beman, W. and Smith, D., editors.

\bibitem[Sho67]{Sho67}
J.R. Shoenfield.
\newblock {\em Mathematical Logic}.
\newblock Addison-Wesley Publishing Company, 1967.

\bibitem[Szm83]{szmielew1983affine}
W.~Szmielew.
\newblock {\em From affine to Euclidean geometry, an axiomatic approach}.
\newblock Polish Scientific Publishers (Warszawa-Poland) and D. Reidel
  Publishing Company (Dordrecht-Holland), 1983.

\bibitem[vS57]{von1857beitrage}
K.G.C. von Staudt.
\newblock {\em Beitr{\"a}ge zur Geometrie der Lage}, volume~2.
\newblock F. Korn, 1857.

\bibitem[Wu94]{bk:Wu1994}
W.-T. Wu.
\newblock {\em Mechanical Theorem Proving in Geometries, Springer 1994}.
\newblock Springer, 1994.
\newblock (Original in Chinese, 1984).

\bibitem[YY05]{YY}
G.C. Young and W.H. Young.
\newblock {\em The first book of geometry}.
\newblock Chelsea Publishing Company, New York, 1905.

\bibitem[Zie82]{ar:ziegler}
Martin Ziegler.
\newblock Einige unentscheidbare {K}\"{o}rpertheorien.
\newblock In V.~Strassen E.~Engeler, H.~L\"{a}uchli, editor, {\em Logic and
  Algorithmic, An international Symposium held in honour of E. Specker}, pages
  381--392. {L'enseignement} {math\'{e}matique}, 1982.

\end{thebibliography}
